\documentclass[12pt]{amsart}
\usepackage[latin1]{inputenc}
\usepackage{amscd}
\usepackage{latexsym}
\usepackage{pstcol,pst-plot,pst-3d}
\usepackage{multicol}
\usepackage{amsmath}
\usepackage{amssymb}
\usepackage{amsthm}
\psset{unit=0.7cm,linewidth=0.8pt,arrowsize=2.5pt 4}

\newpsstyle{fatline}{linewidth=1.5pt}
\newpsstyle{fyp}{fillstyle=solid,fillcolor=verylight}
\definecolor{verylight}{gray}{0.97}
\definecolor{light}{gray}{0.9}
\definecolor{medium}{gray}{0.85}

\def\NZQ{\Bbb}

\def\QQ{{\NZQ Q}}
\def\ZZ{{\NZQ Z}}

\def\frk{\frak}               

\def\mm{{\frk m}}

\def\Phi{{\frk n}}
\def\Phi{{\frk N}}

\def\opn#1#2{\def#1{\operatorname{#2}}} 
\opn\chara{char} \opn\length{\ell} \opn\pd{pd} \opn\rk{rk}
\opn\projdim{proj\,dim} \opn\injdim{inj\,dim} \opn\rank{rank}
\opn\depth{depth} \opn\grade{grade} \opn\height{height}
\opn\embdim{emb\,dim} \opn\codim{codim}

\opn\Tr{Tr} \opn\bigrank{big\,rank}
\opn\superheight{superheight}\opn\lcm{lcm}
\opn\trdeg{tr\,deg}
\opn\reg{reg} \opn\lreg{lreg} \opn\ini{in} \opn\lpd{lpd}
\opn\size{size} \opn\Pf{Pf} \opn\GL{GL} \opn\SL{SL} \opn\mod{mod}
\opn\ord{ord} \opn\Gin{Gin}
\opn\Hilb{Hilb}\opn\adeg{adeg}\opn\std{std}\opn\ip{infpt}
\opn\pol{pol}
\opn\div{div} \opn\Div{Div} \opn\cl{cl} \opn\Cl{Cl}
\opn\Spec{Spec} \opn\Supp{Supp} \opn\supp{supp} \opn\Sing{Sing}
\opn\Ass{Ass} \opn\Min{Min}
\opn\Ann{Ann} \opn\Rad{Rad} \opn\Soc{Soc}
\opn\Syz{Syz} \opn\Im{Im} \opn\Ker{Ker} \opn\Coker{Coker}
\opn\Am{Am} \opn\Hom{Hom} \opn\Tor{Tor} \opn\Ext{Ext}
\opn\End{End} \opn\Aut{Aut} \opn\id{id}

\opn\nat{nat}
\opn\pff{pf}
\opn\Pf{Pf} \opn\GL{GL} \opn\SL{SL} \opn\mod{mod} \opn\ord{ord}
\opn\Gin{Gin} \opn\Hilb{Hilb} \opn\cd{cd}
\opn\aff{aff} \opn\con{conv} \opn\relint{relint} \opn\st{st}
\opn\lk{lk} \opn\cn{cn} \opn\core{core} \opn\vol{vol}
\opn\link{link} \opn\star{star}\opn\limdepth{limdepth}
\opn\limdim{limdim}
\opn\gr{gr}

\def\poly#1#2#3{#1[#2_1,\dots,#2_{#3}]}
\def\pot#1#2{#1[\kern-0.28ex[#2]\kern-0.28ex]}

\opn\dirlim{\underrightarrow{\lim}}
\opn\inivlim{\underleftarrow{\lim}}
\let\iso=\cong

\let\Dirsum=\bigoplus

\let\To=\longrightarrow
\def\Implies{\ifmmode\Longrightarrow \else
        \unskip${}\Longrightarrow{}$\ignorespaces\fi}
\def\implies{\ifmmode\Rightarrow \else
        \unskip${}\Rightarrow{}$\ignorespaces\fi}
\def\iff{\ifmmode\Longleftrightarrow \else
        \unskip${}\Longleftrightarrow{}$\ignorespaces\fi}

\let\:=\colon
\newtheorem{Theorem}{Theorem}[section]
\newtheorem{Lemma}[Theorem]{Lemma}
\newtheorem{Corollary}[Theorem]{Corollary}
\newtheorem{Proposition}[Theorem]{Proposition}

\newtheorem{Definition}[Theorem]{Definition}

\def\Young#1{\vbox{\smallskip\offinterlineskip
    \halign{&\vbox{##}\kern-\Thickness\cr #1}}}

\newdimen\Squaresize \Squaresize=12pt
\newdimen\Thickness \Thickness=.3pt
\newdimen\Correction \Correction=7pt
\let\epsilon\varepsilon
\let\phi=\varphi
\let\kappa=\varkappa
\opn\dis{dis}
\def\pnt{{\raise0.5mm\hbox{\large\bf.}}}
\def\lpnt{{\hbox{\large\bf.}}}
\opn\Lex{Lex}

\def\Coh#1#2{H_{\mm}^{#1}(#2)}

\begin{document}
\title{Tameness of Local cohomology of monomial ideals with respect to monomial prime ideals}
\author{Ahad Rahimi}
\address{Ahad Rahimi, Fachbereich Mathematik und
Informatik, Universit\"at Duisburg-Essen, Campus Essen, 45117
Essen, Germany} \email{ahad.rahimi@uni-essen.de}
\maketitle
\begin{abstract}
In this paper we consider the local cohomology of monomial ideals with respect to monomial prime ideals and show that all these local cohomology modules are tame.
\end{abstract}
\section*{Introduction}
Let $R$ be a graded ring. Recall that a graded $R$-module $N$ is
tame, if there exists an integer $j_0$ such that $N_j=0$ for all
$j\leq j_0$, or else $N_j\neq 0$ for all $j\leq j_0$. Brodmann and
Hellus \cite{BrH} raised the question  whether for a finitely
generated, positively graded algebra $R$ with $R_0$  Noetherian,
 the local cohomology modules $H_{R_+}^i(M)$ for a finitely generated graded $R$-module $M$ are all tame. Here $R_+=\Dirsum_{i>0}R_i$ is  graded irrelevant ideal of $R$. See
\cite{B} for a survey on this problem.

In this paper we only consider rings defined by  monomial
relations. We first consider the squarefree case, since from a
combinatorial point of view this is the more interesting case and
also since in this case the formula which we obtain are more
simple. So let $K$ be a field and $\Delta$ a simplicial complex on
$V=\{v_1,\ldots, v_n\}$. In Section 1, as a generalization of
Hochster's formula \cite{BH}, we compute (Theorem \ref{Hochster})
the Hilbert series of local cohomology of the Stanley-Reisner ring
$K[\Delta]$ with respect to a monomial prime ideal. With the
choice of the monomial prime ideal, the Stanley-Reisner ring
$K[\Delta]$, and hence also all the local cohomology of it, can be
given a natural bigraded structure. In Proposition \ref{dim} we give
a formula for the $K$-dimension of the bigraded components of the
local cohomology modules. Using  this formula we deduce that the local cohomology of
$K[\Delta]$ with respect to a monomial prime ideal is always tame.

 In \cite{T} Takayama  generalized  Hochster's formula to any graded monomial ideal which is not necessarily squarefree. In Section 2, as a
 generalization of Takayama's result,  we compute the Hilbert series of local cohomology of monomial ideals with respect to
 monomial prime ideals and observe  that again all these modules are
 tame. The result proved here is surprising because in a recent
 paper, Cutkosky and Herzog  \cite{CH} gave an example which shows that in general not all
 local cohomology modules are tame.

\section{Local cohomology of Stanley-Reisner rings with respect to monomial prime ideals}
Let $K$ be a field and let $S = \poly{K}{Y}{r}$ be a polynomial
ring with the standard grading.  For a  squarefree monomial ideal
$I\subset S$ we set $R =S/I$.  We denote by $y_i$ the residue
classes of indeterminates $Y_i$ in $R$ for $i=1,\ldots, n$. Thus
we have $R=K[y_1,\dots,y_r]$. We may view $R$ as the
Stanley-Reisner ring of some simplicial complex $\Delta$ with
vertices $\{w_1,\dots,w_r\}$.

Let $P$ be any monomial prime ideal of $R$.  We may assume that
$P=(y_1,\ldots, y_n)$ for some integer $n\leq r$. After this
choice of $P$ we view $R$ as a bigraded $K$-algebra. We rename
some of the variables, and  set  $x_i=y_{n+i}$ for $i=1,\ldots,m$
where $m=r-n$, and assign the following bidegrees:  $\deg
x_i=(1,0)$ for $i=1, \dots, m$ and $\deg y_j=(0,1)$ for
$j=1,\dots, n$.  We decompose the vertex set of  the corresponding
simplicial complex $\Delta$ accordingly, so that $\Delta$ has
vertices $\{v_1,\dots,v_m,w_1,\dots,w_n\}$ where vertices
$V=\{v_1,\dots,v_m\}$ and $W=\{w_1,\dots,w_n\}$ correspond to the
variables of $x_1,\ldots, x_m$ and $y_1,\ldots, y_n$,
respectively. By [1, Theorem 5.1.19] we have
\[
H^i_{P}(R)\iso H^i(C^{\lpnt}) \quad \text {for all} \quad i \geq 0,
\]
where $C^{\lpnt}$ is the \v{C}ech complex

\[
C^{\lpnt} :
0\rightarrow C^0 \rightarrow C^1 \rightarrow \cdots \rightarrow C^n \rightarrow 0
\]
with
\[
C^t = \Dirsum_{1\leq j_1<\cdots <j_t \leq n}R_{y_{j_1}\dots y_{j_t}},
\]
and whose  differential is composed of the maps
\[
(-1)^{s-1} nat : R_{y_{j_1}\dots y_{j_t}} \To R_{y_{j_1}\dots y_{j_{t+1}}},
\]
if $\{i_1,\ldots, i_t\} = \{j_1,\ldots, \hat{j}_s, \ldots, j_{t+1}\}$
and $0$ otherwise.
Note that $C^{\lpnt}$ is a $\ZZ^m \times \ZZ^n$-bigraded complex.
For $(a,b)\in \ZZ^m \times \ZZ^n$ and $y=y_{j_1}\dots y_{j_s}$ with $1\leq j_1<\cdots <j_s \leq n$
 one defines a $\ZZ^m \times \ZZ^n $-bigrading on $R_y$ by setting
\begin{eqnarray}\label{1}
(R_{y})_{(a,b)}=\{r/y^l:  \deg r-l\deg y=(a,b)\}.
\end{eqnarray}
Here $r$ is a bihomogeneous element in $R$, $l\in \ZZ$  and $\deg$
denotes the multi-bidegree. Given
$F=\{w_{j_1},\dots,w_{j_s}\}\subseteq W$ and $b\in \ZZ^{n}$. We
set $G_b=\{w_j: w_1 \leq w_j\leq w_n , b_j<0\}$ ,  $H_b=\{w_j:w_1
\leq w_j\leq w_n , b_j>0\}$  and the support of  $b$ is the set
$\supp b=\{w_j: w_1 \leq w_j\leq w_n , b_j\neq 0\}$. Note that
$\supp b = G_b \cup H_b$.

We set $N_a=\{v_i:v_1 \leq v_i\leq v_m, a_i\neq 0\}=\supp a$ for $a\in \ZZ^m$ and denote by $\ZZ_+^m$ and $\ZZ_-^n$ the sets of
$\{a\in \ZZ^m: a_i\geq 0 \quad \text {for} \quad  i=1,\dots, m\}$ and $\{b\in \ZZ^n: b_i\leq 0 \quad \text {for} \quad  i=1,\dots, n\}$, respectively.
With the notation introduced  one has

\begin{Lemma}
\label{low}
The following statements hold:
\begin{enumerate}
\item[(a)] $\dim_K(R_y)_{(a,b)}\leq 1 $,  for all $a\in \ZZ^{m}$ and $b\in \ZZ^{n}$.
\item[(b)] $(R_y)_{(a,b)}\iso K $,  if and only if $F \supset G_b$, $F\cup H_b\cup N_a\in \Delta$ and $a\in \ZZ_+^m$.
\end{enumerate}
\end{Lemma}
\begin{proof}
As explained before, we may view the standard graded polynomial
ring $S$  as a standard bigraded polynomial ring and then 
$R=K[x_1, \dots , x_m, y_1, \dots, y_n]$ with $m+n=r$ is also
naturally bigraded. Thus part (a) follows from \cite[Lemma 5.3.6
(a)]{BH}. For the proof (b) we set $c=(a,b)$. By \cite[Lemma 5.3.6
(b)]{BH} we have $F \supset G_{c}$ and $F \cup H_{c}\in \Delta$.
Thus (\ref {1}) implies that $a\in \ZZ^m_+$ and hence $G_{c}=G_b$.
We also note that $H_{c}=H_b \cup N_a$.
\end{proof}

As a consequence of Lemma \ref{low}   for $a\in \ZZ^m_+$ , $b\in \ZZ^n$ and $i\in \ZZ$ we observe that $(C^i)_{(a,b)}$ has the following    $K$- basis:
\[
\{b_F: F\supset G_b, F\cup H_b \cup N_a \in \Delta,\left|F\right|=i\}.
\]

Therefore, since $C^{\lpnt}$ is $\ZZ^m \times \ZZ^n$-bigraded
complex one obtains for each $(a,b)\in \ZZ^m \times \ZZ^n$ a
complex
\[
(C^{\lpnt})_{(a,b)} :
0\rightarrow (C^0)_{(a,b)} \rightarrow (C^1)_{(a,b)} \rightarrow \cdots \rightarrow (C^n)_{(a,b)} \rightarrow 0,
\]
of finite dimensional $K$-vector spaces
\[
(C^i)_{(a,b)} =\Dirsum_{{F\supset G_b \atop F\cup H_b \cup N_a \in \Delta}\atop \left|F\right|=i} K b_F.
\]

The differential $\partial:(C^i)_{(a,b)}\longrightarrow (C^{i+1})_{(a,b)}$ is given by $\partial (b_F)=\sum (-1)^{\partial(F,F')}b_{F'}$ where the sum is taken over all $F'$ such that $F'\supset F$,
$F'\cup H_b \cup N_a \in \Delta$ and $\left|F'\right|=i+1$, and where $\partial(F,F')=s$ for $F'=[w_0,\dots,w_i]$ and $F=[w_0,\dots,\hat{w_s},\dots, w_i]$. Then we describe the $(a,b)$th component of the local cohomology in terms of this subcomplex:
\begin{eqnarray}\label{2}
H_{P}^i(K[\Delta])_{(a,b)}=H^i(C^{\lpnt})_{(a,b)} = H^i(C_{(a,b)}^{\lpnt}).
\end{eqnarray}
Let $\Delta$ be a simplicial complex with vertex set $V$ and $\tilde{\mathcal C}(\Delta)$ the augumented oriented chain complex of $\Delta$, see \cite[Section 5.3]{BH}  for details. For an abelian group $G$, the $i$th reduced simplicial cohomology of $\Delta$ with values in $G$ is defined to be  
\begin{eqnarray}\label{3}
\widetilde{H}^i(\Delta;G)=H^i(\Hom_{\ZZ}(\tilde{\mathcal C}(\Delta),G))\quad \text{for all}\quad i.
\end{eqnarray}
 Given  $F\subseteq V$,  we recall the following definitions : The {\em star} of $F$ is the set $st_{\Delta} F=\{ G \in \Delta: F \cup G \in \Delta \}$, and the {\em link} of $F$ is the set $\lk_{\Delta}F=\{G:F \cup G \in \Delta, F \cap G=\emptyset\}$.  We write $st$ and $lk$ instead of $st_\Delta$ and $lk_\Delta$ (for short). We see that  $\st F$ is a subcomplex of $\Delta$, $\lk F$ a subcomplex of $\st F$, and that $\st F=\lk F=\emptyset$ if $F \notin\Delta$.
For $W\subseteq V$, we denote by $\Delta_W$ the simplicial complex restricted to $W$. i.e. the simplicial complex consisting of all faces $F \in \Delta $ whose vertices belong to $W$.

Now in order to compute $H^i(C_{(a,b)}^{\lpnt})$, we prove the following
\begin{Lemma}
\label{link}
For all $a\in \ZZ_+^m$ and $b\in \ZZ^n$ there exists an isomorphism of complexes
\[
(C^{\lpnt})_{(a,b)}\longrightarrow \Hom_{\ZZ}(\tilde{\mathcal C}(\lk_{\st H_b}G_b \cup N_a)_W[-j-1];K),\quad j=\left|G_b\right|
\]
\end{Lemma}
\begin{proof}
The assignment $F\mapsto F'=F-G_b$ establishes a bijection between the set
\[
\beta =\{F\in \Delta_W : F\supset G_b, F\cup H_b \cup N_a \in \Delta,\left|F\right|=i\}.
\]
and the set $\beta' =\{F'\in \Delta_W : F'\in (\lk_{\st H_b } G_b \cup N_a)_W, \left|F'\right|= i-j\}$. Here $F'\in (\lk_{\st H_b } G_b \cup N_a)_W$, since $F'\cap (G_b\cup N_a)=\emptyset$ and $F'\cup (G_b\cup N_a)\in \st H_b$. Therefore  we see that
\[
\alpha ^i:(C^i)_{(a,b)}\longrightarrow  \Hom_{\ZZ}(\tilde{\mathcal C}(\lk_{\st H_b}G_b \cup N_a)_{i-j-1};K),\quad b_F\mapsto \phi_{F-G_b}
\]
is an isomorphism of vector spaces. Here $\phi_{F'}$ is defined by
\[
\phi_{F'}(F'')=\left\{
\begin{array}{cc}
1 & \text{if $F=F''$},\\
0 & \text{otherwise}.
\end{array}
\right.
\]
\end{proof}
As a generalization of Hochster's formula \cite[Theorem 5.3.8]{BH} we prove the following
\begin{Theorem}
\label{Hochster} Let $I \subset S=K[X_1,\ldots, X_m,Y_1,\ldots,
Y_n]$ be a squarefree monomial ideal with the natural $\ZZ^m
\times \ZZ^n$-bigrading. Then the bigraded Hilbert series of the
local cohomology modules of $R=S/I=K[\Delta]$ with respect to the
$\ZZ^m \times \ZZ^n$-bigrading is given by
\[
H_{H_{P}^i(K[\Delta])}(\bold{s},\bold{t})=\sum_{F\in \Delta_W}\sum_{G\subset V}\dim _K \widetilde{H}_{i-\left|F\right|-1}((\lk F \cup G)_W;K)\prod _{v_i\in G}\frac{s_i}{1-s_i} \prod_{w_j\in F}\frac{t_j^{-1}}{1-t_j^{-1}}
\]
where $\bold{s}=(s_1, \dots, s_m)$,  $\bold{t}=(t_1, \dots, t_n)$,  $P=(y_1,\dots,y_n)$ and $\Delta$ is the simplicial complex corresponding to the Stanley-Reisner ring $K[\Delta]$.
\end{Theorem}
\begin{proof}
By (\ref{2}), Lemma \ref{link} and (\ref{3}) we observe that there are isomorphisms of bigraded
$K$-vector spaces
\begin{eqnarray*}
H_{P}^i(K[\Delta])_{(a,b)}& \iso & H^i(\Hom_{\ZZ}(\tilde{\mathcal C}((\lk_{\st H_b}G_b \cup N_a)_W)[-j-1];K),\quad j=\left|G_b\right|\\ & = & H^{i-\left|G_b\right|-1}(\Hom_{\ZZ}(\tilde{\mathcal C}((\lk_{\st H_b}G_b \cup N_a)_W);K)\\ & = &  \widetilde{H}^{i-\left|G_b\right|-1}((\lk_{st H_b} G_b \cup N_a )_W;K),
\end{eqnarray*}
and therefore by \cite[Exercise 5.3.11]{BH} we have
\begin{eqnarray}\label{4}
\dim_K H_{P}^i(K[\Delta])_{(a,b)}=\dim_K \widetilde{H}_{i-\left|G_b\right|-1}((\lk_{st H_b} G_b \cup N_a)_W;K).
\end{eqnarray}
If $H_b \neq \emptyset$  by \cite[Lemma 5.3.5]{BH},  $\lk_{st H_b}G_b \cup N_a$  is acyclic,  and so     $\widetilde{H}_{i-\left|G_b\right|-1}((\lk_{st H_b} G_b \cup N_a)_W;K)=0$  for all $i$.  If $H_b=\emptyset$, then  $st H_b=\Delta$, and so $\lk_{st H_b} G_b \cup N_a=\lk G_b \cup N_a$.
Thus in this case $\supp(b)=G_b$. We also note that  $H_b=\emptyset$ if and only if $b\in \ZZ_-^n$.   In order to simplify notation we will write $s(a)$ and $s(b)$ for the support of $a\in \ZZ^m_+$ and $b\in \ZZ_-^n$, respectively  and set
$d(i, s(b),s(a))=\dim_K \widetilde{H}_{i-\left|s(b)\right|-1}((\lk_{st H_b} s(b) \cup s(a))_W;K)$.
Using these facts and (\ref{4}) we have
\begin{eqnarray*}
H_{H_{P}^i(K[\Delta])}(\bold{s},\bold{t})& = & \sum_{{a\in \ZZ^m_+ , b\in \ZZ_-^n}}\dim_K H_{P}^i(K[\Delta])_{(a,b)}\bold{s}^a \bold{t}^b \\&=&
\sum_{b\in \ZZ^n_-}(\sum_{a\in \ZZ_+^m}d(i,s(b),s(a))\bold{s}^a )\bold{t}^b \\ & =&
\sum_{b\in \ZZ^n_-}(\sum_{G \subset V}\sum_{{s(a)=G} \atop {a\in \ZZ_+^m}}d(i,s(b),s(a))\bold{s}^a)\bold{t}^b \\&=&
\sum_{b\in \ZZ^n_-}(\sum_{G\subset V} d(i,s(b),s(a))\sum_{s(a)=G} \bold{s}^a )\bold{t}^b\\ & = &
\sum_{b\in \ZZ^n_-}(\sum_{G\subset V} d(i,s(b),G )\prod_{v_i \in G}\frac{s_i}{1-s_i})\bold{t}^b\\&=&
 \sum_{F\in \Delta_W}\sum_{{s(b)=F}\atop {b \in \ZZ^n_-}}(\sum_{G\subset V} d(i,s(b),G)\prod_{v_i \in G}\frac{s_i}{1-s_i}) \bold{t}^b \\ & = &
 \sum_{F\in \Delta_W}\sum_{G\subset V}d(i,F,G) \prod _{v_i\in G}\frac{s_i}{1-s_i} \prod_{w_j\in F}\frac{t_j^{-1}}{1-t_j^{-1}},
\end{eqnarray*}
as desired.  Here $\bold s^a=s_1^{a_1}\dots s_m^{a_m}$ for $a=(a_1,\dots,a_m)$ and  $\bold t^b=t_1^{b_1}\dots t_n^{b_n}$ for $b=(b_1,\dots,b_n)$.  We also used the fact that $\sum_{s(a)=G} \bold{s}^a=1$ for $G = \emptyset$  and $ \sum_{s(a)=G} \bold{s}^a=\prod_{v_i \in G}\frac{s_i}{1-s_i}$  for $G\neq \emptyset$.
\end{proof}
We observe that Hochster's formula \cite[Theorem 5.3.8]{BH}
easily follows from Theorem \ref {Hochster}. In fact, if we assume
that $m=0$, then $G=\emptyset$, $(\lk F\cup G)_W=\lk F$ and $\prod
_{v_i\in G}s_i/(1-s_i)=1$. Moreover, we may consider $\deg Y_j=1$
for all $j$. Therefore we get the Hochster formula.

In view of Theorem \ref{Hochster} and (\ref{4}) we get the following isomorphism of $K$-vector spaces
\begin{Corollary}
For all $a\in \ZZ^m_+$ and $b\in \ZZ_-^n$ we have
\[
H_{P}^i(K[\Delta])_{(a,b)}\iso \widetilde{H}^{i-\left|F\right|-1}((\lk F \cup G)_W;K),
\]
where $F=\supp b$ and $G=\supp a$.
\end{Corollary}
\begin{Corollary}
\label{ordinary}
With the notation of Theorem {\em \ref{Hochster}} one has
\begin{eqnarray*}
H_{H_{P}^i(K[\Delta])}(\bold{s},\bold{t})=H_{H_{\mm}^i(K[\Delta_W])}(\bold {t})+ \sum_{F\in \Delta_W}\sum_{G\subset V\atop G\neq \emptyset}d(i,F,G)\prod _{v_i\in G}\frac{s_i}{1-s_i} \prod_{w_j\in F}\frac{t_j^{-1}}{1-t_j^{-1}},
\end{eqnarray*}
where $H_{H_{\mm}^i(K[\Delta_W])}(\bold t)$ is the Hilbert series of  $i$th  ordinary local cohomology of $K[\Delta_W]$ with respect to the maximal ideal $\mm=(y_1,\dots, y_n)$ and where
\[
d(i,F,G)=\dim _K \widetilde{H}_{i-\left|F\right|-1}((\lk F \cup G)_W;K).
\]

\end{Corollary}
\begin{proof}
By Theorem \ref{Hochster} we may write
\begin{eqnarray*}
H_{H_{P}^i(K[\Delta])}(\bold{s},\bold{t})= \sum_{F\in \Delta_W} d(i,F,\emptyset) \prod_{w_j\in F}\frac{t_j^{-1}}{1-t_j^{-1}} &+&\\
 \sum_{F\in \Delta_W}\sum_{G\subset V\atop G\neq \emptyset}d(i,F,G)\prod _{v_i\in G}\frac{s_i}{1-s_i} \prod_{w_j\in F}\frac{t_j^{-1}}{1-t_j^{-1}}.
\end{eqnarray*}
Since
\[
d(i,F,\emptyset) = \dim _K \widetilde{H}_{i-\left|F\right|-1}((\lk
F \cup \emptyset)_W;K)=\dim _K
\widetilde{H}_{i-\left|F\right|-1}((\lk_{\Delta_W} F ;K),
\]
Hochster's formula ( \cite[Theorem  5.3.8]{BH} ) implies the
desired equality.
\end{proof}
In view of Corollary \ref{ordinary} we immediately obtain
\begin{Corollary}
$H_{P}^i(K[\Delta])\neq 0$ for $i=\depth K[\Delta_W] $ and $i=\dim K[\Delta_W]$.
\end{Corollary}

We are interested in the Hilbert series of $H_{P}^i(K[\Delta])$
as a  $\ZZ \times \ZZ$-bigraded algebra. Note that for all $k , j
\in \ZZ$ we have
\begin{eqnarray}\label{5}
H_{P}^i(K[\Delta])_{(k,j)}=\Dirsum_{{a \in \ZZ^m,\left|a\right|=k}\atop {b \in \ZZ^n,\left|b\right|=j}} H_{P}^i(K[\Delta])_{(a,b)},
\end{eqnarray}
where $\left|a\right|=\sum_{i=1}^m a_i$ for $a=(a_1,\dots,a_m)$ and $\left|b\right|=\sum_{i=1}^n b_i$ for $b=(b_1,\dots,b_n)$.

Using this observation we obtain

\begin{Proposition}
\label{dim}
 For all $i$ and $k , j \in \ZZ$  one has
\[
\dim_K H_{P}^i(K[\Delta])_{(k,j)}=\sum_{{F\in \Delta_W}\atop {G\subset V}}d(i,F,G)\binom{k-1}{\left|G\right|-1}\binom{-j-1}{\left|F\right|-1},
\]
where
\[
d(i,F,G)=\dim _K \widetilde{H}_{i-\left|F\right|-1}((\lk F \cup G)_W;K).
\]
\end{Proposition}
\begin{proof}
We set $\left|G\right|=g$ and  $\left|F\right|=f$.
In view of  (\ref{5}) it follows that the Hilbert series of  $H_{P}^i(K[\Delta])$ with respect to the $\ZZ \times \ZZ$-bigraded is obtained from  Theorem \ref {Hochster} by replacing all $s_i$ and $t_j$ by $s$ and $t$, respectively. Thus we have
\[
H_{H_{P}^i(K[\Delta])}(s,t) =  \sum_{F\in \Delta_W} \sum_{G\subset V}d(i,F,G)(\frac{s}{1-s})^g(\frac{t^{-1}}{1-t^{-1}})^f.
\]
We note that
\begin{eqnarray}\label{6}
 \frac{1}{(1-s)^i}=\sum_{r=0}^\infty \binom {i+r-1}{i-1}s^r  \quad \text{for all}\quad  i>0.
\end{eqnarray}
Expanding $(\frac{s}{1-s})^g$ for $g=0$ and $g=1$ and comparing
coefficients with (\ref{6}) we are forced to make the following
convention:  $\binom{-1}{-1}=1$, $\binom{i}{-1}=0$ for all $i\geq
0$ and $\binom{i}{0}=1$ for all $i\geq 0$. Thus we have
\begin{eqnarray*}
 H_{H_{P}^i(K[\Delta])}(s,t) =
\sum_{F\in \Delta_W}\sum_{G\subset V}d(i,F,G)\sum_{r=0}^\infty \binom {g+r-1}{g-1}s^{r+g}\sum_{h=0}^\infty \binom {f+h -1}{f-1}t^{-f-h }.
\end{eqnarray*}
We set $k=r+g$ and $j=-f-h$. Thus $r=k-g$ and $h=-j-f$. Therefore for all $k$ and $j$ with $g\leq k$ and  $0 \leq f \leq -j$ the desired formula follows.
\end{proof}
\begin{Corollary}
For all $i$ and $ j \in \ZZ$ one has
\[
\dim_K H_{P}^i(K[\Delta])_{(0,j)}=\dim_K H_{\mm}^i(K[\Delta_W])_j,
\]
where $H_{\mm}^i(K[\Delta_W])$ is the $i$th ordinary local cohomology of $K[\Delta_W]$ with respect to $\mm=(y_1,\dots,y_n).$
\end{Corollary}
\begin{proof}
By Proposition \ref {dim} and the fact that $(\lk F)_W=\lk_{\Delta_W} F$ we have
\begin{eqnarray*}
\dim_K H_{P}^i(K[\Delta])_{(0,j)} &=& \sum_{F\in \Delta_W}\dim_K \widetilde{H}_{i-\left|F\right|-1}((\lk F)_W;K)\binom{-j-1}{\left|F\right|-1}\\ &=& \sum_{F\in \Delta_W}\dim_K \widetilde{H}_{i-\left|F\right|-1}(\lk_{\Delta_W} F;K)\binom{-j-1}{\left|F\right|-1}\\&=& \dim_K H_{\mm}^i(K[\Delta_W])_j.
\end{eqnarray*}
The last equality follows from Hochster's theorem.
\end{proof}
For all  $j \in \ZZ$, we set \[
H_{P}^i(K[\Delta])_j=\Dirsum_k H_{P}^i(K[\Delta])_{(k,j)},
\]
and consider $H_{P}^i(K[\Delta])_j$ as a finitely generated
graded $R_0$-module. In the following we show that the
Krull-dimension of $H_{P}^i(K[\Delta])_j$ is constant for $j\ll
0$.
\begin{Theorem}
\label{dimtame} For all $i$ there exists an integer $j_0$  such
that for $j\leq j_0$, the Krull-dimension $\dim
H_{P}^i(K[\Delta])_j$ is constant. \end{Theorem}
\begin{proof}

By Proposition \ref {dim} and using (\ref{6}) the $\ZZ$-graded
Hilbert series of $H_{P}^i(K[\Delta])_j$ is given by
\begin{eqnarray*}
H_{H_{P}^i(K[\Delta])_j}(s) & = & \sum_{k=0}^\infty \dim_K H_{P}^i(K[\Delta])_{(k,j)}s^k \\ & = &
 \sum_{k=0}^\infty  \sum_{{F\in \Delta_W}\atop {G\subset V}}d(i,F,G)\binom{k-1}{\left|G\right|-1}\binom{-j-1}{\left|F\right|-1}s^k\\  &=&
  \sum_{{F\in \Delta_W}\atop {G\subset V, \left|G\right|=0}}d(i,F,G)\binom{-j-1}{\left|F\right|-1} \sum_{r=-1}^\infty \binom{r}{-1} s^{r+1}\\&+&
   \sum_{{F\in \Delta_W}\atop {G\subset V, \left|G\right|=1}}d(i,F,G)\binom{-j-1}{\left|F\right|-1} \sum_{r=0}^\infty \binom{r}{0}s^{r+1} \\&+\dots+&
\sum_{{F\in \Delta_W}\atop {G\subset V, \left|G\right|=m}}d(i,F,G)\binom{-j-1}{\left|F\right|-1} \sum_{i=m-1}^\infty \binom{r}{m-1}s^{r+1}\\&=&
A_0(j)+\frac{A_1(j) s}{1-s} + \frac{A_2(j) s^2}{(1-s)^2}+ \dots + \frac{A_m(j) s^m}{(1-s)^m}\\&=&
\frac{\sum_{r=0}^m A_r(j)(1-s)^{m-r}s^r}{(1-s)^m},
\end{eqnarray*}
where
\[
A_r(j)=\sum_{{F\in \Delta_W}\atop {G\subset V, \left|G\right|=r}}d(i,F,G)\binom{-j-1}{\left|F\right|-1}.
\]
We may write
\[
\sum_{r=0}^m A_r(j)(1-s)^{m-r}s^r=\sum_{r=0}^m B_r(j)s^r,
\]
 where  $B_r(j)$ is a polynomial with coefficients in $\QQ$ of degree at most $n-1$.
Therefore we have
\begin{eqnarray}\label{7}
H_{H_{P}^i(K[\Delta])_j}(s)=
\frac{Q_j(s)}{(1-s)^m},
\end{eqnarray}
where $Q_j(s)=\sum_{r=0}^m B_r(j)s^r$. We denote by $Q_j(s)^{(k)}$
the $k$th derivative of $Q_j(s)$ as a function in  $s$ and set
$R_k(j)=[Q_j(s)^{(k)}](1)$ which is of course a polynomial in $j$.
Here we distinguish two cases: First suppose that $R_k=0$ for all
$k\geq 0$. Then the Taylor expansion of $Q_j(s)$
\[
Q_j(s)=\frac{R_0(j)}{0!}+ \frac{R_1(j)}{1!}(1-s)+\frac{R_2(j)}{2!}(1-s)^2+\dots
\]
implies that $Q_j(s)=0$ for all $j$. Thus we see that $R_k=0$ for
all $k\geq 0$, which is equivalent to say that  $Q_j(s)=0$  for
all $j$, and which in turn implies that $H_{P}^i(K[\Delta])_j=0$
for all $j$, and we set  $\dim H_{P}^i(K[\Delta])_j=-\infty$.
Now we assume that not all $R_k=0$, and define
\[
c=\min \{i: R_i\neq 0\}.
\]
 Thus $R_k(j)=0$ for all $j$ and all $k<c$. Since $R_c$ has only finitely many zeroes, it follows that  $R_c(j) \neq 0$ for $j\ll 0$, i.e. there exists an integer $j_0$ such that $R_c(j)\neq 0$ for $j\leq j_0$.
Therefore  $R_k(j)= 0$ for $j\leq j_0$, if $k<c$ and $R_k(j)\neq 0$ for $j\leq j_0$, if $k=c$. Thus  for $j\leq j_0$ we may  write $Q_j(s)=(1-s)^c\tilde{Q}_j(s)$
where $\tilde{Q}_j(s)$ is a polynomial in $s$ with $\tilde{Q}_j(1)\neq 0$. Therefore by (\ref{7}) and \cite[Corollary  4.1.8]{BH}
we have $\dim  H_{P}^i(K[\Delta])_j=m-c$ for all $j\leq j_0$, as desired.
\end{proof}
\begin{Definition}
\label{deftame} {\em Let $R$ be a positively graded Noetherian
ring.  A graded $R$-module $N$ is called {\em tame}, if there exists an integer $j_0$ such that either
\[
N_j=0\quad\text{for all}\quad j\leq j_0,\quad \text{or}\quad N_j\neq 0\quad \text{for all}\quad j\leq j_0.
\]}
\end{Definition}

For example, any finitely generated $R$-module is tame.
\begin{Corollary}
Let $I\subset S=K[ Y_1,\dots, Y_r]$ be a graded squarefree
monomial ideal and let $\Delta$ be the simplicial complex such
that  $K[\Delta]=S/I$. Let $P$ be a monomial prime ideal of
$K[\Delta]$. Then for all $i$, the local cohomology modules of
$H_{P}^i(K[\Delta])$ are  tame.
\end{Corollary}
\begin{proof}
The assertion follows from Theorem \ref{dimtame}.
\end{proof}

\section{Local cohomology of monomial ideals with respect to monomial prime ideals}
 We recall  two results  due to Takayama \cite{T}.
Let $S = \poly{K}{Y}{n}$ be a polynomial ring with the standard
grading.  For a monomial  ideal $I\subset S$ we set $R = S/I$.  We
denote by $y_i$ the image of $Y_i$ in $R$ for $i=1,\ldots, n$ and
set $\mm =(y_1,\ldots, y_n)$, the unique maximal ideal. For
a monomial ideal $I\subset S$, we denote by $G(I)$ the minimal set
of monomial generators. Let $u = Y_1^{c_1}\cdots Y_n^{c_n}$ be a
monomial with $c_j\geq 0$ for all $j$, then we define $\nu_j(u) =
c_j$ for $j=1,\ldots, n$, and $\supp(u)= \{j:  c_j\ne 0\}$. 

We set
$G_b = \{j : b_j < 0\}$ for $b\in\ZZ^n$.  By Takayama we have

\begin{Lemma}
\label{Tak:lemma1}
Let $y= y_{i_1}\cdots y_{i_r}$ with $i_1 < \cdots <i_r$ and set $F = \{i_1,\ldots, i_r\}$.
For all $b\in\ZZ^n$ we have $\dim_K (R_y)_b \leq 1$ and the following are equivalent
\begin{enumerate}
\item [$(i)$]  $(R_y)_b \iso K$
\item [$(ii)$]  $F \supset G_b$ and for all $u\in G(I)$ there exists $j\notin F$ such that $\nu_j(u) > b_j\geq 0.$
\end{enumerate}
\end{Lemma}
For any $b\in\ZZ^n$, we define a simplicial complex
\begin{eqnarray*}
\Delta_b = \left\{ F- G_b
    \;\vert\;
    \begin{array}{l}
     F \supset G_b, \mbox{ and for all $u\in G(I)$ there exists $j\notin F$}\\
     \mbox{ such that } \nu_j(u) > b_j\geq 0
                \end{array}
\right\}.
\end{eqnarray*}
\begin{Theorem}
\label{Tak}
Let $I\subset S = \poly{K}{Y}{n}$ be a monomial ideal. Then
the multigraded Hilbert series of the local cohomology modules of $R = S/I$
with respect to the $\ZZ^n$-grading is given by
\begin{equation*}
\Hilb(\Coh{i}{R}, {\bf t})
=\sum_{F\in\Delta}
  \sum
   \dim_K\tilde{H}_{i-\vert F\vert -1}(\Delta_b; K) {\bf t}^b
\end{equation*}
where ${\bf t} = (t_1,\cdots, t_n)$. The second sum runs over
$a\in\ZZ^n$ such that $G_b = F$ and $b_j\leq \rho_j-1$,
$j=1,\ldots,n$, with $\rho_j = \max\{\nu_j(u):  u\in G(I)\}$ for
$j=1,\ldots, n$, and $\Delta$ is the simplicial complex
corresponding to the Stanley-Reisner ideal $\sqrt{I}$.
\end{Theorem}
Note that Takayama's formula can be rewritten as following
\begin{eqnarray*}
\Hilb(\Coh{i}{R}, {\bf t}) & = & \sum_{F\in\Delta}
\dim_K\tilde{H}_{i-\vert F\vert -1}(\Delta_b; K) \sum_{{b \in \ZZ^m,  G_b=F \atop  b_j\leq \rho_j-1} \atop j=1,\dots, m } {\bf t}^b \\ & =& \sum_{F\in\Delta}
\dim_K\tilde{H}_{i-\vert F\vert -1}(\Delta_b; K)\prod_{w_j\in F} \frac{1-{t_j}^{\rho_j}}{1-t_j}\prod_{w_j\notin F}\frac{{t_j}^{-1}}{1-{t_j}^{-1}}.
\end{eqnarray*}
We see that in the squarefree case, this formula together with
\cite[Corollary 1]{T} implies again Hochster's formula.

Now let $S = K[X_1,\dots, X_m, Y_1,\dots, Y_n]$  be a standard
bigraded polynomial ring over $K$.  For a monomial  ideal
$I\subset S$ we set $R =S/I$.  The residue classes of the
variables will be denoted by  $x_i$ and $y_j$  and set $P
=(y_1,\ldots, y_n)$. For monomial $u\in S$ we may write $u=u_1u_2$
where $u_1$ and $u_2$ are monomials in $X$ and $Y$.  Let $\Delta$
be the simplicial complex corresponding to $\sqrt{I}$.  As before
we denote the vertices corresponding to the $X_i$ by $v_i$ and
those corresponding to the $Y_j$ by $w_j$. We set $G_b = \{w_j :
b_j < 0\}$ for $b\in\ZZ^n$. With the same arguments as Lemma
\ref{Tak:lemma1} we have
\begin{Lemma}
\label{Tak2}
Let $y= y_{i_1}\cdots y_{i_r}$ with $i_1 < \cdots <i_r$ and set $F =\{i_1,\dots,i_r\}$.
For all $a\in\ZZ^m$ and $b\in\ZZ^n$ we have $\dim_K (R_y)_{(a,b)} \leq 1$ and the following are equivalent
\begin{enumerate}
\item [$(i)$]  $(R_y)_{(a,b)} \iso K$ \item [$(ii)$]  $F \supset
G_b$,  $a\in \ZZ^m_+$,  and for all $u\in G(I)$ there exists
$j\notin F$ such that $\nu_j(u_2) > b_j\geq 0$ or  for at least
one $i$, $\nu_i(u_1) > a_i\geq 0$.
\end{enumerate}
\end{Lemma}
For any $a\in\ZZ^m_+$ and  $b\in\ZZ^n$, we define a simplicial
complex
\begin{eqnarray*}
\Delta_{(a,b)} = \left\{ F- G_b
    \;\vert\;
    \begin{array}{l}
     F \supset G_b,  a\in\ZZ^m_+ \mbox{ and for all $u\in G(I)$ there exists $j\notin F$} \\ \mbox{ such that }  \nu_j(u_2) > b_j\geq 0,     \mbox{or for at least one i}, \nu_i(u_1) > a_i\geq 0
                \end{array}
\right\}.
\end{eqnarray*}
As a generalization of Takayama's result we have
\begin{Theorem}
\label{Hochster 2} Let $I\subset S=K[X_1,\dots, X_m, Y_1,\dots,
Y_n]$ be a monomial ideal with the natural $\ZZ^m
\times \ZZ^n$-bigrading. Then the bigraded Hilbert series
of the local cohomology modules of $R = S/I$ with respect to the
$\ZZ^m \times \ZZ^n$-bigrading is given by

\begin{eqnarray*}
 & &  \Hilb(H^i_{P}(R){\bf, s,t}) =  \\   &  &
     \sum_{F\in\Delta_W}\sum_{G\subset V}
D(i,F,G) \prod_{v_i\notin G} \frac{1-{s_i}^{-\sigma_i}}{1-(s_i)^{-1}}\prod_{v_i\in G}\frac{s_i}{1-s_i}\prod_{w_j\notin F} \frac{1-{t_j}^{\rho_j}}{1-t_j}\prod_{w_j\in F}\frac{{t_j}^{-1}}{1-{t_j}^{-1}},
\end{eqnarray*}
where
$P=(y_1,\dots,y_n)$, $
D(i,F,G)=\dim_K\tilde{H}_{i-\vert F\vert -1}(\Delta_{(a,b)}; K),
$
${\bf s} = (s_1,\cdots, s_m)$, ${\bf t} =( t_1, \cdots, t_n)$,
$\rho_j = \max\{\nu_j(u_2)\: u\in G(I)\}$ for $j=1,\ldots, n$,
$\sigma_i = \max\{\nu_i(u_1): u\in G(I)\}$
for $j=1,\ldots, m$, and $\Delta$ is the
simplicial complex corresponding to the Stanley-Reisner
ideal $\sqrt{I}$.
\end{Theorem}
\begin{proof}
With the same arguments as  in the proof  Theorem 1 in \cite{T} we
can show that
\begin{equation*}
\Hilb(H_{P}^i(R), {\bf s , t})
=\sum \sum \dim_K\tilde{H}_{i-\vert F\vert -1}(\Delta_{(a,b)}; K) {\bf s}^a {\bf t}^b
\end{equation*}
where the first sum runs over $F\in \Delta_W$, $b\in\ZZ^n$ such
that $G_b = F$ and $b_j\leq \rho_j-1$, $j=1,\ldots,n$, and the
second sum runs over $a\in\ZZ^m$ such that $N_a = G$ and $a_i\geq
\sigma_i-1$, $i=1,\ldots,m$. Indeed, if  we assume that $b_j >
\rho_j-1$ and $a_i < \sigma_i-1$, the proof of \cite[Theorem
1]{T} shows that $\dim_K\tilde{H}_{i-\vert F\vert
-1}(\Delta_{(a,b)}; K)=0$ for all $i$. Therefore we may write
\begin{eqnarray*}
& &\Hilb(H^i_{P}(R){\bf, s,t})  =  \sum_{F\in\Delta_W}\sum_{G\subset V}
D(i,F,G) \sum_{{a \in \ZZ^m,  N_a=G \atop  a_i\geq \sigma_i-1} \atop i=1,\dots, m } {\bf s}^a\sum_{{b \in \ZZ^n,
 G_b=F \atop  b_j\leq \rho_j-1} \atop j=1,\dots, n } {\bf t}^b,
\end{eqnarray*}
where
\[
D(i,F,G)=\dim_K\tilde{H}_{i-\vert F\vert -1}(\Delta_{(a,b)}; K).
\]
Since
\[
\sum_{{a \in \ZZ^m,  N_a=G \atop  a_i\geq \sigma_i-1} \atop i=1,\dots, m } {\bf s}^a=\prod_{v_i\notin G} \frac{1-{s_i}^{-\sigma_i}}{1-(s_i)^{-1}}\prod_{v_i\in G}\frac{s_i}{1-s_i}
\]
and
\[
\sum_{{b \in \ZZ^n,  G_b=F \atop  b_j\leq \rho_j-1} \atop
j=1,\dots, n } {\bf t}^b=\prod_{w_j\notin F}
\frac{1-{t_j}^{\rho_j}}{1-t_j}\prod_{w_j\in
F}\frac{{t_j}^{-1}}{1-{t_j}^{-1}},
\]
the desired formula follows.
\end{proof}
We observe that Theorem \ref{Hochster} is a special case of Theorem \ref{Hochster 2}. In fact, if we assume that
$\sigma_i=1$ for $i=1,\dots, m$ and $\rho_j=1$ for $j=1,\dots, n$, then $a\in \ZZ^m_+$, $b\in \ZZ^n_-$,
$\prod_{v_i\notin G} \frac{1-{s_i}^{-\sigma_i}}{1-(s_i)^{-1}}=1$ and $\prod_{w_j\notin F} \frac{1-{t_j}^{\rho_j}}{1-t_j}=1$. Moreover, by the proof of \cite[Corollary  1]{T} we have
$\Delta_{(a,b)}=\lk_{\st H_{(a,b)}} G_{(a,b)}=\lk_{\st N_a\cup H_b}G_b=\lk_{\st N_a}G_b=(\lk F \cup G)_W$.
\begin{Proposition}
\label{dim 2}
 For all $i$ and $k , j \in \ZZ$  one has
\[
\dim_K H_{P}^i(R)_{(k,j)}=\sum_{{F\in \Delta_W}\atop {G\subset V}}D(i,F,G)\sum_{r=0}^ {\sigma_G} a_G(r)\binom{k+r-1}{\left|G\right|-1}\sum_{h=0}^ {\rho_F} b_F(h)\binom{h-j-1}{\left|F\right|-1},
\]
where
$
\sigma_G=\sum_{v_i\notin G}(\sigma_i-1)$, $\rho_F=\sum_{w_j\notin F}(\rho_j-1)$ and  $a_G(r), b_F(h)\in \ZZ$ for $r=0,\dots, \sigma_G$ and $h=0, \dots, \rho_F$.
\end{Proposition}
\begin{proof}
In Theorem \ref{Hochster 2} we replace all $s_i$ by $s$ and all
$t_j$ by $t$, and obtain
\[
H_{H_{P}^i(R)}(s,t) =  \sum_{F\in \Delta_W} \sum_{G\subset V}D(i,F,G)P_G(s^{-1})(\frac{s}{1-s})^{\left|G\right|}Q_F(t)(\frac{t^{-1}}{1-t^{-1}})^{\left|F\right|},
\]
where
\[
P_G(s^{-1})=\prod_{v_i\notin G}(1+s^{-1}+\dots+ s^{-\sigma_i+1})\quad and \quad Q_F(t)=\prod_{w_j\notin F}(1+t^{1}+\dots+ t^{\rho_j-1})
\]
with $\deg P_G(s^{-1})=\sigma_G$ and
$\deg Q_F(t)=\rho_F.$
We may write $P_G(s^{-1})=\sum_{r=0}^{\sigma_G} a_G(r) s^{-r}$ where $a_G(r)\in \ZZ$ for $r=0,\dots, \sigma_G$ and $Q_F(t)=\sum_{h=0}^{\rho_F} b_F(h) t^{h}$ where $b_F(h)\in \ZZ$ for $h=0, \dots, \rho_F$.
By setting $\left|G\right|=g$ and $\left|F\right|=f$,  we have

\[ H_{H_{P}^i(R)}(s,t) =
 \sum_{F\in \Delta_W} \sum_{G\subset
V}D(i,F,G)A_G(s)B_F(t),
\]
where
\[A_G(s)=\sum_{r=0}^{\sigma_G} a_G(r) \sum_{l=0}^\infty \binom
{g+l-1}{g-1}s^{l+g-r}
\]
and
\[
B_F(t)=\sum_{h=0}^{\rho_F} b_F(h)\sum_{\rho=0}^\infty \binom
{f+\rho -1}{f-1}t^{-f-\rho +h}.
\]

We set $k=l+g-r$ and $j=-f-\rho +h$. Then  $r=l+g-k$ and
$h=j+f+\rho $, and the desired formula follows.
\end{proof}

\begin{Theorem}
\label{dimtame2}
For all $i$ there exists an integer $j_0$  such that for $j\leq j_0$, the Krull- dimension $\dim H_{P}^i(R)_j$ is constant. \end{Theorem}
\begin{proof}

By Proposition \ref {dim 2} the $\ZZ$-graded Hilbert series of $H_{P}^i(R)_j$ is given by
\begin{eqnarray*}
H_{H_{P}^i(R)_j}(s) & = & \sum_{k=0}^\infty \dim_K H_{P}^i(R)_{(k,j)}s^k \\ & = &
 \sum_{k=0}^\infty  \sum_{{F\in \Delta_W}\atop {G\subset V}}D(i,F,G)\sum_{t=0}^{\sigma_G}a_G(t) \binom{k+t-1}{\left|G\right|-1}\sum_{h=0}^{\rho_F}b_F(h) \binom{h-j-1}{\left|F\right|-1}s^k\\  &=&
   A_0(j)\sum_{r=-1}^\infty \sum_{t=0}^{\sigma_G}a_G(t)\binom{t+r}{-1} s^{r+1}  \\&+&
   A_1(j) \sum_{r=-1}^\infty \sum_{t=0}^{\sigma_G}a_G(t)\binom{t+r}{0}s^{r+1} \\&+\dots+&
A_m(j)\sum_{r=-1}^\infty
\sum_{t=0}^{\sigma_G}a_G(t)\binom{t+r}{m-1}s^{r+1}\\&=&
A_0(j)+A_1(j)\frac{P_1(s)}{1-s} + A_2(j)\frac{P_2(s)}{(1-s)^2}+
\dots + A_m(j)\frac{P_m(s)}{(1-s)^m}\\&=& \frac{\sum_{r=0}^m
A_r(j)(1-s)^{m-r}P_r(s)}{(1-s)^m},
\end{eqnarray*}
where
\[
A_r(j)=\sum_{{F\in \Delta_W}\atop {G\subset V,
\left|G\right|=r}}D(i,F,G)\sum_{h=0}^{\rho_F}b_F(h)\binom{h-j-1}{\left|F\right|-1},  
\]
and $P_r(s)\in  \ZZ[s]$ with $\deg P_r(s)= r.$ Here we used (\ref{6}) and that $a_G(0)=1$ for $G=\emptyset$,
$\binom{t+r}{-1}=1$ for $t+r=-1$ and $0$ otherwise, $\binom{t+r}{0}=1$ for $t+r\geq0$ and
$ \binom{t+r}{n}=\binom{t}{n}+r\binom{t}{n-1}+ \frac{r(r-1)}{2!}\binom{t}{n-2}+\frac{r(r-1)(r-2)}{3!}\binom{t}{n-3}+\dots+ \frac{r(r-1)(r-2)\dots (r-n+1)}{n!}\binom{t}{0}$  for $t+r \geq n \geq 1.$
We may write
\[
\sum_{r=0}^m A_r(j)(1-s)^{m-r}P_r(s)=\sum_{r=0}^m B_r(j)s^r,
\]
 where  $B_r(j)$ is a polynomial with coefficients in $\QQ$ of degree at most $n-1$.
Therefore we have
\[
H_{H_{P}^i(R)_j}(s)=
\frac{Q_j(s)}{(1-s)^m},
\]
where $Q_j(s)=\sum_{r=0}^m B_r(j)s^r$. We proceed in the same way
as in the proof of Theorem \ref{dimtame} and get the desired
result.
\end{proof}
\begin{Corollary}
Let $I\subset S=K[Y_1,\dots, Y_r]$ be any monomial ideal and set $R=S/I.$ Let $P$ be a monomial prime ideal in $R$. Then for all $i$, the local cohomology modules of $H_{P}^i(R)$  are tame.
\end{Corollary}
\begin{proof}
The assertion follows from Theorem \ref {dimtame2}.
\end{proof}

\begin{center}
{\bf Acknowledgment} 
\end{center}
\hspace*{\parindent}
I am grateful to Professor Jürgen Herzog for numerous helpful comments and discussions

\end{document}